\documentclass[11pt]{amsart}
\usepackage{amssymb,latexsym,amsmath,amscd,fullpage,epsfig,mathrsfs,yfonts}
\input xy
\xyoption{all}

\newcommand{\RR}{\mathbb R}
\newcommand{\ZZ}{\mathbb Z}
\newcommand{\QQ}{\mathbb Q}
\newcommand{\CC}{\mathbb C}

\newcommand{\N}{\mathbb N}

\def\Hom{\textup{Hom}}
\def\rk{\textup{rk}}
\def\dim{\textup{dim}}

\def\qgr{\textup{qgr}}
\def\colim{\textup{colim}}
\def\Proj{\textup{Proj}}
\def\proje{\textup{proj}}
\def\Tor{\textup{Tor}}
\def\Torsion{\textup{Torsion}}
\def\QCoh{\textup{QCoh}}
\def\Coh{\textup{Coh}}
\def\Mod{\textup{Mod}}
\def\cohproj{\textup{cohproj}}
\def\deg{\textup{deg}}
\def\GKdim{\textup{GKdim}}
\def\Ext{\textup{Ext}}
\def\End{\textup{End}}
\def\Lim{\textup{Lim}}
\def\Gr{\textup{Gr}}
\def\QGr{\textup{QGr}}

\newcommand{\proj}{\mathbb P}

\newcommand{\spc}{\vspace{3mm}}
\newcommand{\Sheaf}{\mathscr}
\newcommand{\beq}{\begin{eqnarray}}
\newcommand{\beqn}{\begin{eqnarray*}}
\newcommand{\eeq}{\end{eqnarray}}
\newcommand{\eeqn}{\end{eqnarray*}}
\newcommand{\nc}{non-commutative \hspace{0.5mm}}

\newcommand{\QED}{\begin{flushright}$\blacksquare$\end{flushright}}
\newcommand{\ul}{\underline}

\newtheorem{thm}{Theorem}[section]
\newtheorem{fact}{Fact}

\newtheorem{lem}[thm]{Lemma}
\newtheorem{prop}[thm]{Proposition}

\newtheorem{ex}{Example}
\newtheorem{defn}[thm]{Definition}
\theoremstyle{remark}
\newtheorem{rem}[thm]{Remark}

\newcommand{\sheaf}{\mathcal{O}}
\newcommand{\map}{\longrightarrow}

\newcommand{\lk}{{}_{A}k}

\begin{document}

%----------------------------------------------------------------------
\title{Lecture notes on non-commutative algebraic geometry and noncommutative tori}
\author{Snigdhayan Mahanta}
\address{Max-Planck-Institut f\"ur Mathematik, Bonn, Germany.}

%----------------------------------------------------------------------

%\begin{abstract}
%In these set of notes I have outlined my understanding of some parts
%of {\it A. Polishchuk's} effort in \cite{Pol2} to reconcile the
%operator algebraic approach towards \nc geometry with that of {\it
%Artin} and {\it Zhang} via module categories. I have restricted my
%attention to the more algebraic aspects of Polishchuk's article,
%which should be complemented by the talks of {\it Jorge Plazas} and
%{\it Behrang Noohi}. The first section is dedicated to the {\it
%Artin} and {\it Zhang} approach followed by a rather hasty (due to
%lack of time) overview of \cite{Pol2}.
%\end{abstract}

%----------------------------------------------------------------------
%\keywords{asd}%

%\date{}%
%\dedicatory{}%
%\commby{}%
%----------------------------------------------------------------------

\maketitle

%--------------------------today---------------------------------------
\begin{center}
\today
\end{center}

%------------------it is time to begin---------------------------------

%======================================================================

\spc
\begin{center}
{\bf INTRODUCTION}
\end{center}

\spc

I would like to thank all the organizers, namely, {\it M. Khalkhali, M. Marcolli,
M. Shahshahani} and {\it M. M. Sheikh-Jabbari}, of the {\it International
Workshop on Noncommutative Geometry, 2005} for giving me the opportunity to speak.

\vspace{5mm}

In section 1 we shall browse through some interesting definitions and
constructions which will be referred to later on. In section 2 we shall discuss \nc projective geometry as initiated by {\it Artin} and {\it Zhang} in
\cite{AZ}. This section is rather long and the readers can easily skip over
some details. In section 3 we shall provide a brief overview of the algebraic
aspects of noncommutative tori, which comprise the most widely studied class
of noncommutative differentiable manifolds. These notes are not
entirely self-contained and should be read in tandem with those of {\it
  B. Noohi} for background material and of {\it J. Plazas} for a better
understanding of the topological and differential aspects of \nc tori.

%-----------------After Artin and Zhang-----------------------------------

\newpage

\section{Some Preliminaries}

\vspace{5mm}

In a paper entitled {\it Some Algebras Associated to Automorphisms of Elliptic
  Curves} Artin, Tate and Van den Bergh \cite{ATV1} gave a nice description of
  \nc algebras which should, in principle, be algebras of functions of some
  nonsingular ``\nc schemes''. In the commutative case, nonsingularity is
  reflected in the regularity of the ring. However, this notion is
  insufficient for \nc purposes. So Artin and Schelter gave a stronger
  regularity condition which we call the {\it Artin-Schelter}(AS)-{\it regularity} condition. The main result of the above-mentioned paper says that AS-regular algebras of dimension 3 (global dimension) can be described neatly as some algebras associated to automorphisms of projective schemes, mainly elliptic curves. Also such algebras are both left and right Noetherian. This subsection is entirely based on the contents of \cite{ATV1}.
\vspace{2mm}

To begin with, we fix an algebraically closed field $k$ of
characteristic $0$. We shall mostly be concerned with $\N$-graded
$k$-algebras $A\;=\;\underset{i\geqslant 0}{\oplus} A_i$, that are
finitely generated in degree $1$, with $A_0$ finite dimensional as a
$k$-vector space. Such algebras are called {\it finitely graded} for
short, though the term could be a bit misleading at first sight.
A finitely graded ring is called {\it connected graded} if
$A_0\;=\;k$. $A_{+}$ stands for the two-sided {\it augmentation
ideal} $\underset{i>0}{\oplus} A_i$.

\spc

\begin{defn} (AS-regular algebra)

A connected graded ring $A$ is called Artin-Schelter (AS) regular of dimension $d$ if it satisfies the following conditions:

(1) $A$ has global dimension $d$.

(2) \text{\GKdim}($A$) $<$ $\infty$.

(3) $A$ is AS-Gorenstein.

\end{defn}

It is worthwhile to say a few words about {\it Gelfand-Kirillov\;dimension} (GKdim) and the {\it AS-Gorenstein} condition of algebras.

\spc

Take any connected graded $k$-algebra $A$ and choose a finite dimensional $k$-vector space $V$ such that $k[V]\;=\;A$. Now set $F^{n}A\;=\;k\;+\;\sum_{i=1}^{n} V^i$ for $n \geqslant 1$. This defines a filtration of $A$. Then the $\GKdim(A)$ is defined to be

\begin{eqnarray*}
\GKdim(A)\;=\;\underset{n}{\textup{lim sup}}\frac{\textup{ln}(\textup{dim}_k F^{n}A)}{\textup{ln}(n)}.
\end{eqnarray*}

Of course, one has to check that the definition does not depend on the choice of $V$.

\begin{rem}
{\it Bergman} \cite{KL} has shown that $\GKdim$ can be any real number $\alpha \geqslant 2$. However, if $\GKdim\;\leqslant\;2$, then it is either $0$ or $1$.
\end{rem}

There are some equivalent formulations of the  AS-Gorenstein condition available in literature. We shall be content by saying the following:

\newpage

\begin{defn} (AS-Gorenstein condition) \label{Gorenstein}

A connected graded $k$-algebra $A$ of global dimension $d<\infty$ is AS-Gorenstein if
\begin{eqnarray*}
\Ext^i_A(k,A)\;=\;0\;for\;i\neq d\;and\;\Ext^d_A(k,A)\simeq k
\end{eqnarray*}
\end{defn}

All regular commutative rings are AS-Gorenstein, which supports our conviction that the AS-Gorenstein hypothesis is desirable for \nc analogues of regular commutative rings. Further, note that the usual {\it Gorenstein} condition (for commutative rings) requires that they be Noetherian of finite injective dimension as modules over themselves.

Now we take up the task of describing the minimal projective resolution $(\star)$ of an AS-regular algebra of dimension $d\;=\;3$. As a fact, let us also mention that the global dimension of a graded algebra is equal to the projective dimension of the left module $\lk$. Let

\beq
0 \map P^d\map\cdots\overset{f_2}{\map}P^1\overset{f_1}{\map}P^0\map \lk\map 0
\eeq

be a minimal projective resolution of the left module $\lk$. $P^0$ turns out
to be $A$; $P^1$ and $P^2$ need an investigation into the structure of $A$ for their descriptions. Suppose $A\;=\;T/I$, where $T\;=\;k\{x_1,\dots,x_n\}$ is a free associative algebra generated by homogeneous elements $x_i$ with degrees $l_{1j}$ (also assume that $\{x_1,\dots,x_n\}$ is a minimal set of generators). Then

\beq
P^1\;\approx\;\underset{j=1}{\overset{n}{\oplus}}\;A(-l_{1j})
\eeq

The map $P^1\map P^0$, denoted $x$, is given by right multiplication with the column vector $(x_1,\dots,x_n)^t$.

Coming to $P^2$, let $\{f_j\}$ be a minimal set of homogeneous generators for the ideal $I$ such that $\deg\;f_j\;=\;l_{2j}$. In $T$, write each $f_j$ as

\beq
f_j\;=\;\underset{j}{\sum}m_{ij}x_j
\eeq

where $m_{ij}\;\in\;T_{l_{2i}-l_{1j}}$. Let $M$ be the image in $A$ of the matrix $(m_{ij})$. Then

\beq
P^2\;\approx\;\underset{j}{\oplus} A(-l_{2j})
\eeq

and the map $P^2\map P^1$, denoted $M$, is just right multiplication by the matrix $M$.

\spc

In general, it is not so easy to interpret all the terms of the resolution $(1)$. However, for a regular algebra of dimension $3$, the resolution looks like

\beq
0 \map A(-s-1)\stackrel{x^t}{\map} A(-s)^r\stackrel{M}{\map} A(-1)^r\stackrel{x}{\map} A\map \lk \map 0 \hspace{10mm}
\eeq

where $(r,s)\;=\;(3,2)$ or $(2,3)$.  Thus such an algebra has $r$ generators and $r$ relations each of degree $s$, $r+s\;=\;5$. Set $g\;=\;(x^t)M$; then

\beq
g^t\;=\;((x^t)M)^t\;=\;QMx\;=\;Qf \hspace{10mm}
\eeq

for some $Q \in GL_r(k)$.

Now, with some foresight, we introduce a new definition, that of a {\it standard algebra}, in which we extract all the essential properties of AS-regular algebras of dimension $3$.

\begin{defn} An algebra $A$ is called \underline{standard} if it can be presented by $r$ generators $x_j$ of degree $1$ and $r$ relations $f_i$ of degree $s$, such that, with M defined by $(3)$, $(r,s)\;=\;(2,3)$ or $(3,2)$ as above, and there is an element $Q \in GL_r(k)$ such that $(6)$ holds.
\end{defn}

\begin{rem}
For a standard algebra $A$, $(5)$ is just a complex and if it is a resolution, then $A$ is a regular algebra of dimension 3.
\end{rem}

\newpage

\begin{center}
{\bf Twisted Homogeneous Coordinate Rings}
\end{center}

\spc

Here we sketch a very general recipe for manufacturing interesting
\nc rings out of a completely ``commutative geometric'' piece of
datum, called an {\it abstract triple}, which turns out to be an
isomorphism invariant for ``AS-regular algebras".

\begin{defn}
An abstract triple $\mathcal{T}=(X,\sigma,\Sheaf{L})$ is a triple consisting of a projective scheme $X$, an automorphism $\sigma$ of
$X$ and an invertible sheaf $\Sheaf{L}$ on $X$.
\end{defn}

It is time to construct the {\it Twisted Homogeneous Coordinate
Ring} $B(\mathcal{T})$ out of an abstract triple. For each integer
$n\geqslant 1$ set

\beq
\Sheaf{L}_n\;=\;\Sheaf{L}\otimes\Sheaf{L}^\sigma\otimes\dots\otimes\Sheaf{L}^{\sigma^{n-1}}
\eeq

where $\Sheaf{L}^\sigma\; :=\;\sigma^{\ast}\Sheaf{L}$. The tensor
products are taken over ${\sheaf}_X$ and we set
$\Sheaf{L}_0\;=\;{\sheaf}_X$. As a graded vector space,
$B(\mathcal{T})$ is defined as

\begin{eqnarray*}
B(\mathcal{T})\;=\;\underset{n\geqslant 0}{\oplus}
H^0(X,\Sheaf{L}_n)
\end{eqnarray*}

For every pair of integers $m,n\geqslant 0$, there is a canonical
isomorphism

\beqn \Sheaf{L}_m\otimes_{{}_k} \Sheaf{L}_n^{\sigma^m} \map
\Sheaf{L}_{m+n} \eeqn

and hence it defines a multiplication on $B(\mathcal{T})$

\beqn H^0(X,\Sheaf{L}_m)\otimes_{{}_k} H^0(X,\Sheaf{L}_n)\map
H^0(X,\Sheaf{L}_{m+n}). \eeqn

\begin{ex}

Let us compute (more precisely, allude to the computation of) the
twisted homogeneous coordinate ring in a very simple case. Let
$\mathcal{T}=(\proj^1,\sheaf(1),\sigma)$, where $\sigma(a_0,a_1)=(q
a_0,a_1)$ for some $q\in k^\ast$. It is based on our understanding
of Example 3.4 of \cite{StV}.

\end{ex}

We may choose a parameter $u$ for $\proj^1$, so that the standard affine
open cover consisting of $U=\proj^1\setminus \{\infty\}$ and $V=\proj^1\setminus \{0\}$ has rings of regular functions $\sheaf(U)=k[u]$ and
$\sheaf(V)=k[u^{-1}]$ respectively. Now we can identify $\sheaf(1)$ with the
sheaf of functions on $\proj^1$ which have at most a simple pole at
infinity; in other words, it is the subsheaf of $k(u)=k(\proj^1)$
generated by $\{1,u\}$. It can be checked that
$H^0(X,\sheaf(n))$ is spanned by $\{1,u,\dots,u^n\}$ and that, as a
graded vector space $B(\proj^1,\textup{id},\sheaf(1))=k\{x,y\}$ (the free
algebra over $k$ generated by $x$ and $y$ and not the usual
polynomial ring), where $x=1$ and $y=u$, thought of as elements of
$B_1=H^0(X,\sheaf(1))$. It should be mentioned that $\sigma$ acts on
the rational functions on the right as $f^\sigma(p)=f(\sigma(p))$
for any $f\in k(\proj^1)$ and $p\in \proj^1$. From the presentation of the algebra it is evident that $\sheaf(1)^\sigma\cong\sheaf(1)$. So as a
graded vector space $B(\mathcal{T})\cong B(\proj^1,id,\sheaf(1))$.
The multiplication is somewhat different though.

\beqn y\cdot x\; =\; y\otimes x^\sigma\; =\; u\otimes 1^\sigma\; =\;
u\otimes 1\; =\; u\in H^0(\proj^1,\sheaf(2)). \eeqn

On the other hand,

\beqn x\cdot y\; =\; x\otimes y^\sigma\; =\; 1\otimes u^\sigma\; =\;
1\otimes qu\; =\; qu\in H^0(\proj^1,\sheaf(2)). \eeqn

So we find a relation between $x$ and $y$, namely, $x\cdot y
\,-\,qy\cdot x\,=\,0$ and a little bit more work shows that this is
the only relation. So the twisted homogeneous coordinate ring
associated to $\mathcal{T}$ is

\beqn B(\mathcal{T})\;=\; k\{x,y\}/(x\cdot y\,-\,qy\cdot x). \eeqn

\noindent A good reference for a better understanding of these rings
is \cite{AV}.

\newpage

\vspace{3mm}
\begin{center}
{\bf A cursory glance at Grothendieck Categories}
\end{center}

For the convenience of the reader let me say a few things about a {\it generator} of a category. An object $G$ of a category $\mathcal{C}$ is called a {\it generator} if, given a pair of morphisms $f,g:A \map B$ in $\mathcal{C}$ with $f\neq g$, there exists an $h:G\map A$ with $fh\neq gh$ (more briefly, $\Hom(G,{}_{-}):\mathcal{C} \map Set$ is a faithful functor). A family of objects $\{G_i\}_{i \in I}$ is called a {\it generating set} if, given a pair of morphisms $f,g:A{}_{\map}^{\map} B$ with $f \neq g$, there exists an $h_i:G_i \map A$ for some $i \in I$ with $fh_i \neq gh_i$. Strictly speaking, this is a misnomer. In a cocomplete category ({\it i.e.,} closed under all coproducts), a family of objects $\{G_i\}_{i\in I}$ forms a generating set if and only if the coproduct of the family forms a generator.

\begin{rem}
Let $\mathcal{C}$ be a cocomplete abelian category. Then an object $G$ is a generator if and only if, for any object $A \in Ob(\mathcal{C})$, there exists an epimorphism

\begin{eqnarray*}
G^{\oplus I} \map A
\end{eqnarray*}

for some indexing set $I$.

\end{rem}

\begin{defn}
\noindent
\text{\bf Grothendieck Category}

  It is a (locally small) cocomplete abelian category with a generator and satisfying, for every family of short
  exact sequences indexed by a filtered category $I$ [{\it i.e.,} $I$ is
  nonempty and, if $i,j\in Ob(I)$ then $\exists$ $k\in Ob(I)$ and arrows
  $i\map k$ and $j\map k$, and for any two arrows $i\,
  {}^{\overset{u}{\map}}_{\underset{v}{\map}}\, j$ $\exists$ $k\in Ob(I)$ and
  an arrow $w:j\map k$ such that $wu = wv$ (think of a categorical formulation
  of a directed set)].

\begin{eqnarray*}
0 \map A_i \map B_i \map C_i \map 0
\end{eqnarray*}

the following short sequence is also exact

\begin{eqnarray*}
0 \map \underset{\underset{i\in I}{\map}}{\colim}\,A_i \map \underset{\underset{i\in I}{\map}}{\colim}\,B_i \map \underset{\underset{i\in I}{\map}}{\colim}\,C_i \map 0
\end{eqnarray*}

{\it i.e.,}  passing on to filtered colimits preserves exactness. This is equivalent to the sup condition of the famous {\bf AB5 Property}.

\end{defn}

\begin{rem}

The original {\bf AB5 Property} requires the so-called {\it sup} condition, besides cocompleteness. An abelian category satisfies {\it sup} if

for any ascending chain $\Omega$ of subobjects of an object $M$, the supremum of $\Omega$ exists; and for any subobject $N$ of $M$, the canonical morphism

\begin{eqnarray*}
sup\;\{L\cap N | L\in \Omega\} \overset{\sim}{\map} (sup\;\Omega) \cap N
\end{eqnarray*}

is an isomorphism. Hence, another definition of a Grothendieck category could be a cocomplete abelian category, having a generator and satisfying the sup condition, {\it i.e.,} an {\bf AB5} category with a generator.
\end{rem}

\begin{ex}
\textbf{(Grothendieck Categories)}\newline

1. $\Mod(R)$, where $R$ is an associative ring with unity.

2. The category of sheaves of $R$-modules on an arbitrary topological space.

3. In the same vein, the category of abelian pre-sheaves on a {\it site} $\mathcal{T}$. Actually this is just $Funct(\mathcal{T}^{op},Ab)$.

4. $\QCoh(X)$, where $X$ is a quasi-compact and quasi-separated scheme. (a morphism of schemes $f:X\map Y$ is called quasi-compact if, for any open quasi-compact $U\subseteq Y$, $f^{-1}(U)$ is quasi-compact in $X$ and it becomes quasi-separated if the canonical morphism $\delta_f:X\map X\times_Y X$ is quasi-compact. A scheme $X$ is called quasi-compact (resp. quasi-separated) if the canonical unique morphism $X\map Spec(\ZZ)$, $Spec(\ZZ)$ being the final object, is quasi-compact (resp. quasi-separated)).
\end{ex}

Grothendieck categories have some remarkable properties which make them amenable to homological arguments.

1. Grothendieck categories are complete {\it i.e.,} they are closed under products.

2. In a Grothendieck category every object has an injective envelope, in particular there are enough injectives.

\vspace{2mm}

The deepest result about Grothendieck categories is given by the following theorem.

\begin{thm} \text{[Gabriel,Popescu]} \label{Gabriel_Popescu}
Let $\mathcal{C}$ be a Grothendieck category and let $\mathcal{G}$ be a generator of $\mathcal{C}$. Put $S = \End(\mathcal{G})$. Then the functor

\begin{eqnarray*}
\Hom(\mathcal{G},-):\mathcal{C}\map \Mod(S^{op})
\end{eqnarray*}

is fully faithful (and has an exact left adjoint).
\end{thm}

\spc

\begin{center}
{\bf Justification for bringing in Grothendieck Categories}
\end{center}

\vspace{2mm}

We begin by directly quoting {\it Manin} \cite{Man2} - ``...Grothendieck taught us, to do geometry you really don't need a space, all you need is a category of sheaves on this would-be space.'' This idea gets a boost from the following reconstruction theorem.

\begin{thm}\text{[Gabriel,Rosenberg \cite{Ros2}]}  \label{reconstruction}
Any scheme can be reconstructed uniquely up to isomorphism from the category of quasi-coherent sheaves on it.
\end{thm}

\begin{rem}
When it is known in advance that the scheme to be reconstructed is an affine
one, we can just take the centre of the category, which is the endomorphism
ring of the identity functor of the category. More precisely, let $X =
Spec\,A$ be an affine scheme and let $\mathcal{A}$ be the category of
quasi-coherent sheaves on $X$, which is the same as $\Mod(A)$. Then the centre
of $\mathcal{A}$, denoted $\End(Id_\mathcal{A})$, is canonically isomorphic to
$A$. [The centre of an abelian category is manifestly commutative and, in general, it gives us only the centre of the ring, that is, $\mathcal{Z}(A)$. But here we are talking about honest schemes and hence $\mathcal{Z}(A)=A$].
\end{rem}

We can also get a derived analogue of the above result, which is, however, considerably weaker. Also it is claimed to be an easy consequence of the above theorem in \cite{BO1}.

\begin{thm} \text{[Bondal,Orlov \cite{BO1}]}
Let X be a smooth irreducible projective variety with ample canonical or anti-canonical sheaf. If $\mathcal{D} = D^b \Coh(X)$ is equivalent as a graded category to $D^b \Coh(X')$ for some other smooth algebraic variety $X'$, then $X$ is isomorphic to $X'$.
\end{thm}

\begin{rem}
Notice that for elliptic curves or, in general, abelian varieties the theorem above is not applicable.
\end{rem}

Finally, consider a pre-additive category with a single object, say $\ast$. Then, being a pre-additive category, $\Hom(\ast,\ast)$ is endowed with an abelian group structure. If we define a product on it by composition, then it is easy to verify that the two operations satisfy the ring axioms. So $\Hom(\ast,\ast)$ or simply $\End(\ast)$ is a ring and that is all we need to know about the pre-additive category. Extrapolating this line of thought, we say that pre-additive categories generalise the concept of rings and since schemes are concocted from commutative rings, it is reasonable to believe that a pre-additive category with some geometric properties should give rise to a ``\nc scheme''. The geometric properties desirable in an abelian category were written down by {\it Grothendieck} as the famous {\bf AB Properties} in \cite{Gro}.

\newpage

\begin{center}
{\bf A brief discussion on construction of quotient categories}
\end{center}

\vspace{2mm}

This discussion is based on {\it Gabriel's} article \cite{Gab} (page 365) and curious readers are encouraged to go through the details from there.

Recall that we call a full subcategory $\mathcal{C}$ of an abelian category $\mathcal{A}$ {\it thick} if the following condition is satisfied:

\vspace{1mm}

\hspace{4mm} for all short exact sequences in $\mathcal{A}$ of the form $0\map M'\map M\map M''\map 0$, we have

\begin{eqnarray*}
M\in \mathcal{C} \Longleftrightarrow \text{ both } M',M''\in \mathcal{C}
\end{eqnarray*}

\vspace{1mm}

Now we construct the quotient of $\mathcal{A}$ by a thick subcategory $\mathcal{C}$, denoted $\mathcal{A/C}$, as follows:

\begin{eqnarray*}
\begin{split}
Ob(\mathcal{A/C}) &= \text{objects of $\mathcal{A}$.} \\
\Hom_{\mathcal{A/C}}(M,N) &= \underset{\underset{\underset{s.t.\, M/M',N'\in\mathcal{C}}{\forall\,M',N'\,\textup{subobj.}}}{\longrightarrow}}{\Lim}\Hom_{\mathcal{A}}(M',N/N').
\end{split}
\end{eqnarray*}

One needs to check that as $M'$ and $N'$ run through all subobjects
of $M$ and $N$ respectively, such that $M/M'$ and $N'$ are in
$Ob(\mathcal{C})$ (take intersection and sum respectively), the
abelian groups $\Hom_{\mathcal{C}}(M',N/N')$ form a directed system.
It satisfies the obvious universal properties which the readers are
invited to formulate. It comes equipped with a canonical quotient
functor $\pi:\mathcal{A}\map\mathcal{A}/\mathcal{C}$.

\begin{prop} \text{[Gabriel]}
Let $\mathcal{C}$ be a thick subcategory of an abelian category $\mathcal{A}$. Then the category $\mathcal{A}/\mathcal{C}$ is abelian and the canonical functor $\pi$ is exact.
\end{prop}

The essence of this quotient construction is that the objects of $\mathcal{C}$ become isomorphic to zero.

\begin{ex}
Let $\mathcal{A} = \Mod(\ZZ)$ and $\mathcal{C} = \Torsion\;groups$.
Then one can show that $\mathcal{A}/\mathcal{C} \simeq \Mod(\QQ)$.

Let us define a functor from $\mathcal{A}/\mathcal{C}$ to $\Mod(\QQ)$
by tensoring with $\QQ$. We simplify the $\Hom$ sets of
$\mathcal{A}/\mathcal{C}$. Using the structure theorem, write every
abelian group as a direct sum of its torsion part and torsion-free
part. If one of the variables is torsion, it can be shown that in
the limit $\Hom$ becomes $0$. So we may assume that both variables
are torsion-free and for simplicity let us consider both of them to
be $\ZZ$. Then,

\beqn
\begin{split}
\Hom_{\mathcal{A}/\mathcal{C}}(\ZZ,\ZZ)&=\underset{n}{{\underset{\map}{{\Hom}_{\mathcal{A}}}}} (n\ZZ,\ZZ) \\
                             &=\underset{n}{\bigcup}\;\frac{1}{n}\ZZ \\
                             &=\QQ = \Hom(\QQ,\QQ)
\end{split}
\eeqn

This says that the functor is full, and an easy verification shows that it is faithful and essentially surjective.

\end{ex}

\vspace{2mm}

      Now we are ready to discuss a model of \nc projective geometry after
      {\it Artin} and {\it Zhang} \cite{AZ}. We would also like to bring to the notice of readers the works of {\it Verevkin} (see \cite{Ver}). But before that let us go through one nice result in the affine case. Let $X$ be an affine scheme and put $A = \Gamma(X,\sheaf_X)$. Then it is well-known that $\QCoh(X)$ is equivalent to $\Mod(A)$. This fact encourages us to ask: which Grothendieck categories can be written as $\Mod(A)$ for some possibly \nc ring $A$?

The answer to this question is given by the theorem below.

\begin{thm} \text{\cite{Ste}}
Let $\mathcal{C}$ be a Grothendieck category with a projective generator $G$ and assume that $G$ is small [{\it i.e.,} $\Hom(G,-)$ commutes with all direct sums]. Then $\mathcal{C} \simeq \Mod(A^{op})$, for $A = \End(G)$.
\end{thm}

Note that the Gabriel-Popescu Theorem [\ref{Gabriel_Popescu}] gave just a fully faithful embedding with an exact left adjoint and not an equivalence.

\newpage

\section{Non-commutative Projective Geometry}

Fix an algebraically closed field $k$; then we shall mostly be dealing with categories which are {\it $k$-linear} abelian categories [{\it i.e.,}
      the bifunctor $\Hom$ ends up in $\Mod(k)$]. Since in commutative algebraic geometry one mostly deals with finitely generated $k$-algebras, which are Noetherian, here we assume that our $k$-algebras are at least {\it right Noetherian}. Later on we shall need to relax this Noetherian condition, but for now we stick to it. Let $R$ be a graded algebra. Then we introduce some categories here:

\begin{eqnarray*}
\begin{split}
\QCoh(X) :=\; &\text{category of quasi-coherent sheaves on a scheme $X$.} \\
\Coh(X)  :=\; &\text{category of coherent sheaves on $X$.} \\
\Mod(A)  :=\; &\text{category of right $A$-modules, where $A$ is a $k$-algebra.} \\
\Gr(R)   :=\; &\text{category of $\ZZ$-graded right $R$-modules, with degree $0$ morphisms.} \\
\Tor(R)  :=\; &\text{full subcategory of $\Gr(R)$ generated by torsion modules}\\
           & \text{({\it i.e.,} $M$ such that $\forall$ $x\in M$,
  $xR_{\geqslant s} = 0$ for some $s$), which is thick.} \\
\QGr(R)  :=\; &\text{the quotient category $\Gr(R)/\Tor(R)$ (refer to the quotient } \\
             &\text{construction before).}
\end{split}
\end{eqnarray*}

\spc

Notice that $\QCoh(X)$ is not obtained from $\Coh(X)$ by a quotient construction as $\QGr(R)$ is from $\Gr(R)$. In fact, when $X$ is Noetherian, $\Coh(X)$ is the subcategory of $\QCoh(X)$ generated by all Noetherian objects in it.

\spc

\begin{rem} \label{notation}
{\bf Standard Convention}. If $XYuvw(\dots)$ denotes an abelian category, then
we shall denote by $xyuvw(\dots)$ the full subcategory consisting of
Noetherian objects and if $A,B,\dots,M,N,\dots$ denote objects in $\Gr(R)$ then
we shall denote by $\mathcal{A},\mathcal{B},\dots,\mathcal{M},\mathcal{N},\dots$ the corresponding objects in $\QGr(R)$.

Some people denote $\QGr(R)$ by $Tails(R)$, but we shall stick to our notation. We denote the quotient functor $\Gr(R)\map \QGr(R)$ by $\pi$. It has a right adjoint functor $\omega : \QGr(R)\map \Gr(R)$ and so, for all $M\in \Gr(R)$ and $\Sheaf{F}\in \QGr(R)$ one obtains

\beqn
\Hom_{\QGr(R)}(\pi M,\Sheaf{F}) \cong \Hom_{\Gr(R)}(N,\omega\Sheaf{F}).
\eeqn

The $\Hom$'s of $\QGr(R)$ take a more intelligible form with the assumptions on
$R$. It turns out that for any $N\in gr(R)$ and $M\in \Gr(R)$,

\beq \label{Hom_desc}
\Hom_{\QGr(R)}(\pi N,\pi M) \cong \underset{\rightarrow}{\Lim}
\Hom_{\Gr(R)}(N_{\geqslant n},M)
\eeq

For any functor $F$ from a $k$-linear category $\mathcal{C}$ equipped with an autoequivalence $s$, we denote by $\underline{F}$ the graded analogue of $F$ given by $\underline{F}(A):=\underset{n\in\ZZ}{\oplus}F(s^nA)$ for any $A\in Ob(\mathcal{C})$. Further, to simplify notation we shall sometimes denote $s^nA$ by $A[n]$ when there is no chance of confusion.
\end{rem}

\newpage

\noindent
Keeping in mind the notations introduced above we have

\begin{lem}
$\omega\pi M \cong \underset{\rightarrow}{\Lim} \ul{\Hom}_R (R_{\geqslant n}, M)$
\end{lem}

{\bf Sketch of proof:}

\beqn
\begin{split}
\omega\pi M & = \ul{\Hom}_R (R,\omega\pi M) \hspace{3mm} \text{[since $R\in
  gr(R)$]} \\
            & = \underset{d\in \ZZ}{\oplus} \Hom_{\Gr(R)} (R,\omega\pi M[d]) \\
            & = \underset{d\in \ZZ}{\oplus} \Hom_{\QGr(R)}(\pi R,\pi M[d])
  \hspace{3mm} \text{[by adjointness of $\pi$ and $\omega$]} \\
            & = \underset{d\in \ZZ}{\oplus} \underset{\rightarrow}{\Lim}
            \Hom_{\Gr(R)} (R_{\geqslant n},M[d]) \hspace{3mm} \text{[by
            \eqref{Hom_desc}]} \\
            & = \underset{\rightarrow}{\Lim} \underset{d\in \ZZ}{\oplus}
            \Hom_{\Gr(R)} (R_{\geqslant n}, M[d]) \\
            & = \underset{\rightarrow}{\Lim} \ul{\Hom}_R (R_{\geqslant n},M)
\end{split}
\eeqn

\QED

The upshot of this lemma is that there is a natural equivalence of functors $\omega \simeq \ul{\Hom} (\pi R,{}_{-})$.

\vspace{2mm}
\begin{center}
\underline{\bf $\Proj\, R$}
\end{center}

\vspace{3mm}

Let $X$ be a projective scheme with a line bundle $\Sheaf{L}$. Then
the homogeneous coordinate ring $B$ associated to $(X,\Sheaf{L})$ is
defined by the formula $B = \underset{n\in
\mathbb{N}}{\oplus}\Gamma(X,\Sheaf{L}^{n})$ with the obvious
multiplication. Similarly, if $\Sheaf{M}$ is a quasi-coherent sheaf
on $X$, $\Gamma_h(\Sheaf{M})=\underset{n\in
\mathbb{N}}{\oplus}\Gamma(X,\Sheaf{M}\otimes\Sheaf{L}^{n})$ defines
a graded $B$-module. Thus, the compostion of $\Gamma_h$ with the
natural projection from $\Gr(B)$ to $\QGr(B)$ yields a functor
${\bar{\Gamma}}_h:\QCoh(X)\map \QGr(B)$. This functor works
particularly well when $\Sheaf{L}$ is ample, as is evident from the
following fundamental result due to {\it Serre}.

\begin{thm} \cite{Ser}
1. Let $\Sheaf{L}$ be an ample line bundle on a projective scheme $X$. Then the functor ${\bar{\Gamma}}_h(-)$ defines an equivalence of categories between $\QCoh(X)$ and $\QGr(B)$.

2. Conversely, if $R$ is a commutative connected graded $k$-algebra, that is, $R_0 = k$ and it is generated by $R_1$ as an $R_0$-algebra, then there exists a line bundle $\Sheaf{L}$ over $X = \Proj(R)$ such that $R = B(X,\Sheaf{L})$, up to a finite dimensional vector space. Once again, $\QGr(R)\simeq \QCoh(X)$.
\end{thm}

In commutative algebraic geometry one defines the $\Proj$ of a graded
ring to be the set of all homogeneous prime ideals which do not
contain the augmentation ideal. This notion is not practicable over
arbitrary algebras. However, {\it Serre's} theorem filters out the
essential ingredients to define the $\Proj$ of an arbitrary algebra.
The equivalence is controlled by the category $\QCoh(X)$, the
structure sheaf $\sheaf_X$ and the autoequivalence given by
tensoring with $\Sheaf{L}$, which depends on the polarization of
$X$. Borrowing this idea we get to the definition of $\Proj$.
Actually one should have worked with a $\ZZ$-graded algebra $R$ and
defined its $\Proj$ but it has been shown in \cite{AZ} that, with the
definition to be provided below, $\Proj\, R$ is the same as $\Proj\,
R_{\geqslant 0}$. Hence, we assume that $R$ is an
$\mathbb{N}$-graded $k$-algebra. $\Gr(R)$ has a shift operator $s$
such that $s(M)=M[1]$ and a special object, $R_R$. We can actually
recover $R$ from the triple $(\Gr(R),R_R,s)$ by

\begin{eqnarray*}
R = \underset{i\in\mathbb{N}}{\oplus} \Hom(R_R,s^{i}(R_R))
\end{eqnarray*}

and the composition is given as follows: $a\in R_i$ and $b\in R_j$,
then $ab = s^{j}(a)\circ b\in R_{i+j}$.

 Let $\mathcal{R}$ denote the image of $R$ in $\QGr(R)$ and we continue to denote by $s$ the autoequivalence induced by $s$ on $\QGr(R)$.

\begin{defn} ($\Proj\, R$)

The triple $(\QGr(R),\mathcal{R},s)$ is called the projective scheme of $R$
and is denoted $\Proj\,R$. Keeping in mind our convention we denote
$(\qgr(R),\mathcal{R},s)$ by $\proje\,R$. This is just as good because there is a way to switch back and forth between $\QGr(R)$ and $\qgr(R)$.
\end{defn}

\begin{center}
\underline{\bf Characterization of $\Proj\, R$}
\end{center}

\vspace{3mm}

We have simply transformed {\it Serre's} theorem into a definiton. It is time
to address the most natural question: which triples
$(\mathcal{C},\mathcal{A},s)$ are of the form $\Proj(R)$ for some graded
algebra $R$? This problem of characterization has been dealt with
comprehensively by {\it Artin} and {\it Zhang}. We will be content with just taking a quick look at the important points. Let us acquaint ourselves with morphisms of such triples. A morphism between $(\mathcal{C},\mathcal{A},s)$ and $(\mathcal{C}',\mathcal{A}',s')$ is given by a triple $(f,\theta,\mu)$, where $f:\mathcal{C}\map\mathcal{C}'$ is a $k$-linear functor, $\theta:f(\mathcal{A})\map \mathcal{A}'$ is an isomorphism in $\mathcal{C}'$ and $\mu$ is a natural isomorphism of functors $f\circ s\map s'\circ f$. The question of characterization is easier to deal with when $s$ is actually an automorphism of $\mathcal{C}$. To circumvent this problem, an elegant construction has been provided in \cite{AZ} whereby one can pass to a different triple, where $s$ becomes necessarily an automorphism. If $s$ is an automorphism one can take negative powers of $s$ as well and it becomes easier to define the graded analogues of all functors (refer to [\ref{notation}]). Sweeping that discussion under the carpet, henceforth, we tacitly assume that $s$ is an automorphism of $\mathcal{C}$ (even though we may write $s$ as an autoequivalence).

 The definition of $\Proj$ was conjured up from {\it Serre's} theorem where the triple was $(\QCoh(X),\sheaf_X ,{}_{-}\otimes\Sheaf{L})$. Of course, one can easily associate a graded $k$-algebra to any $(\mathcal{C},\mathcal{A},s)$.

\begin{eqnarray} \label{triple2alg}
\Gamma_h(\mathcal{C},\mathcal{A},s) = \underset{n\geqslant 0}{\oplus} \Hom(\mathcal{A},s^n\mathcal{A})
\end{eqnarray}

with multiplication $a\cdot b = s^n(a)b$ for $a\in \Hom(\mathcal{A},s^m\mathcal{A})$ and $b\in \Hom(\mathcal{A},s^n\mathcal{A})$.

\begin{rem}
Let $X$ be a scheme, $\sigma\in Aut(X)$ and $\Sheaf{L}$ be a line bundle on
$X$. Then one obtains the {\it twisted homogeneous coordinate ring}, as
discussed in section 1, as a special case of the above construction applied to
the triple $(\QCoh(X),\sheaf_X,\sigma_\ast(- \otimes\Sheaf{L}))$. [Hint: to verify this, use the projection formula for sheaves]
\end{rem}

Notice that $\Sheaf{L}$ needs to be ample for {\it Serre's} theorem to work. So we need a notion of ampleness in the categorical set-up.

\begin{defn} (Ampleness) \label{ample}

Assume that $\mathcal{C}$ is locally Noetherian. Let $\mathcal{A}\in Ob(\mathcal{C})$ be a Noetherian object and let $s$ be an autoequivalence of $\mathcal{C}$. Then the pair $(\mathcal{A},s)$  is called ample if the following conditions hold:

1. For every Noetherian object $\mathcal{O}\in Ob(\mathcal{C})$ there are positive integers $l_1,\dots,l_p$ and an epimorphism from $\overset{p}{\underset{i=0}{\oplus}}\mathcal{A}(-l_i)$ to $\mathcal{O}$.

2. For every epimorphism between Noetherian objects $\mathcal{P}\map\mathcal{Q}$ the induced map $\Hom(\mathcal{A}(-n),\mathcal{P})\map \Hom(\mathcal{A}(-n),\mathcal{Q})$ is surjective for $n\gg 0$.
\end{defn}

\begin{rem}
The first part of this definition corresponds to the standard definition of an ample sheaf and the second part to the homological one.
\end{rem}

Now we are in good shape to state one part of the theorem of {\it Artin} and {\it Zhang} which generalises that of {\it Serre}.

\begin{thm} \label{characterization}
Let $(\mathcal{C},\mathcal{A},s)$ be a triple as above such that the following conditions hold:

(H1) $\mathcal{A}$ is Noetherian,

(H2) $A := \Hom(\mathcal{A},\mathcal{A})$ is a right Noetherian ring and $\Hom(\mathcal{A},\mathcal{M})$ is a finite $A$-module for all Noetherian $\mathcal{M}$, and

(H3) $(\mathcal{A},s)$ is ample.

Then $\mathcal{C}\simeq \QGr(B)$ for $B=\Gamma_h (\mathcal{C},\mathcal{A},s)$. Besides, $B$ is right Noetherian.
\end{thm}

The converse to this theorem requires an extra hypothesis, which is the so-called $\chi_1$ condition. One could suspect, and rightly so, that there is a $\chi_n$ condition for every $n$. They are all some kind of condition on the graded $\Ext$ functor. However, they all look quite mysterious at first glance. Actually most naturally occurring algebras satisfy them but the reason behind their occurrence is not well understood. We shall discuss them in some cases later but we state a small proposition first.

\begin{prop} \label{chi_equiv}

Let $M\in \Gr(B)$ and fix $i\geqslant 0$. There is a right $B$-module structure on $\underline{\Ext}^n_B (B/B_+,M)$ coming from the right $B$-module structure of $B/B_+$. Then the following are equivalent:

1. for all $j\leqslant i$, $\ul{\Ext}^j_B (B/B_{+},M)$ is a finite $B$-module;

2. for all $j\leqslant i$, $\ul{\Ext}^j_B (B/B_{\geqslant n},M)$ is finite for all $n$;

3. for all $j\leqslant i$ and all $N\in \Gr(B)$, $\ul{\Ext}^j_B (N/N_{\geqslant n},M)$ has a right bound independent of $n$;

4. for all $j\leqslant i$ and all $N\in \Gr(B)$, $\underset{\rightarrow}{\Lim}\ul{\Ext}^j_B (N/N_{\geqslant n},M)$ is right bounded.
\end{prop}

       The proof is a matter of unwinding the definitions of the terms
       suitably and then playing with them. We shall do something smarter
       instead - make a definition out of it.

\begin{defn} ($\chi$ conditions)

A graded algebra $B$ satisfies $\chi_n$ if, for any finitely generated graded $B$-module $M$, one of the equivalent conditions of the above proposition is satisfied (after substituting $i = n$ in them). Moreover, we say that $B$ satisfies $\chi$ if it satisfies $\chi_n$ for every $n$.
\end{defn}

\begin{rem}

Since $B/B_+$ is a finitely generated $B_0$-module ($B_0 = k$) we could have equally well required the finiteness of $\underline{\Ext}^n_B (B/B_+,M)$ over $B_0 = k$ {\it i.e.,} $\dim_k \underline{\Ext}^n_B (B/B_+,M) < \infty$ for $\chi_n$.
\end{rem}

Let $B$ be an $\mathbb{N}$-graded right Noetherian algebra and $\pi : \Gr(B)\map \QGr(B)$.

\begin{thm}
If $B$ satisfies $\chi_1$ as well, then (H1), (H2) and (H3) hold for the triple $(\qgr(B),\pi B,s)$. Moreover, if $A=\Gamma_h (\QGr(B),\pi B,s)$, then $\Proj\,B$ is isomorphic to $\Proj\,A$ via a canonical homomorphism $B\map A$. [We have a canonical map $B_n = \Hom_B(B,B[n])\map \Hom(\pi B,\pi B[n]) = A_n$ given by the functor $\pi$.]
\end{thm}

The proofs of the these theorems are once again quite long and involved. So they are omitted. What we need now is a good cohomology theory for studying such \nc projective schemes.

\begin{center}
\underline{\bf Cohomology of $\Proj\, R$}
\end{center}

\vspace{3mm}

The following rather edifying theorem due to {\it Serre} gives us some insight into the cohomology of projective (commutative) spaces.

\begin{thm} \cite{Har}
Let $X$ be a projective scheme over a Noetherian ring $A$, and let $\sheaf_X(1)$ be a very ample invertible sheaf on $X$ over $Spec\,A$. Let $\Sheaf{F}$ be a coherent sheaf on $X$. Then:

1. for each $i\geqslant 0$, $H^i(X,\Sheaf{F})$ is a finitely generated $A$-module,

2. there is an integer $n_0$, depending on $\Sheaf{F}$, such that for each $i>0$ and each $n\geqslant n_0$, $H^i(X,\Sheaf{F}(n))=0$.
\end{thm}

There is an analogue of the above result and we zero in on that. We have already come across the $\chi$ conditions, which have many desirable consequences. Actually the categorical notion of ampleness doesn't quite suffice. For the desired result to go through, we need the algebra to satisfy $\chi$ too. Without inundating our minds with all the details of $\chi$ we propose to get to the point {\it i.e.,} cohomology. Set $\pi R =\mathcal{R}$. On a projective (commutative) scheme $X$ one can define the sheaf cohomology of $\Sheaf{F}\in \Coh(X)$ as the right derived functor of the global sections functor, {\it i.e.,} $\Gamma$. But $\Gamma(X,\Sheaf{F})\cong \Hom_{\sheaf_X}(\sheaf_X,\Sheaf{F})$. Buoyed by this fact, the proposed definition of the cohomology for every $\mathcal{M}\in \qgr(R)$ is

\beqn
H^n(\mathcal{M}) := \Ext^n_{\mathcal{R}}(\mathcal{R},\mathcal{M})
\eeqn

 However, taking into consideration the graded nature of our objects we also define the following:

\begin{eqnarray*}
\underline{H}^n(\Proj\,R,\mathcal{M}) := \ul{\Ext}^n_{\mathcal{R}}(\mathcal{R},\mathcal{M}) = \underset{i\in\ZZ}{\oplus} \Ext^n_{\mathcal{R}}(\mathcal{R},\mathcal{M}[i])
\end{eqnarray*}

The category $\QGr(R)$ has enough injectives and one can choose a nice ``minimal'' injective resolution of $\mathcal{M}$ to compute its cohomologies, the details of which are available in chapter 7 of \cite{AZ}.

Let $M\in \Gr(R)$ and write $\mathcal{M} = \pi M$. Then one should observe that

\beq    \label{cohom_desc}
\begin{split}
\ul{H}^n(\mathcal{M}) & = \ul{\Ext}^n_{\QGr(R)}(\mathcal{R},\mathcal{M}) \\
                      & \cong \underset{\rightarrow}{\Lim}\ul{\Ext}^n_{\Gr(R)}
                      (R_{\geqslant n},M) \hspace{3mm} [by \eqref{Hom_desc}]
\end{split}
\eeq

\vspace{2mm}
As $R$-modules we have the following exact sequence,

\begin{eqnarray*}
0\map R_{\geqslant n}\map R\map R/R_{\geqslant n}\map 0
\end{eqnarray*}

For any $M\in \Gr(R)$, the associated $\underline{\Ext}$ sequence in $\Gr(R)$ looks like

\begin{eqnarray*}
\dots\underline{\Ext}^j(R/R_{\geqslant n},M)\map\underline{\Ext}^j(R,M)\map\underline{\Ext}^j(R_{\geqslant n},M)\map\dots
\end{eqnarray*}

Since $R$ is projective as an $R$ module, $\underline{\Ext}^j(R,M)=0$ for every
$j\geqslant 1$. Thus, we get the following exact sequence:

\beq \label{exactseq}
\hspace{10mm} 0\rightarrow\underline{\Hom}(R/R_{\geqslant n},M)\rightarrow M\rightarrow\underline{\Hom}(R_{\geqslant n},M)\rightarrow\underline{\Ext}^1(R/R_{\geqslant n},M)\rightarrow 0
\eeq

and, for every $j\geqslant 1$, an isomorphism

\beq \label{Ext_relation}
\underline{\Ext}^j(R_{\geqslant n},M)\cong\underline{\Ext}^{j+1}(R/R_{\geqslant n},M)
\eeq

The following theorem is an apt culmination of all our efforts.

\begin{thm}  (Serre's finiteness theorem)

Let $R$ be a right Noetherian $\mathbb{N}$-graded algebra satisfying $\chi$, and let $\Sheaf{F}\in \qgr(R)$. Then,

(H4) for every $j\geqslant 0$, ${H}^j(\Sheaf{F})$ is a finite right $R_0$-module, and

(H5) for every $j\geqslant 1$, $\ul{H}^j(\Sheaf{F})$ is right bounded; {\it i.e.,} for $d \gg 0$, ${H}^j(\Sheaf{F}[d])=0$.
\end{thm}

\noindent
{\bf Sketch of proof:}

Write $\Sheaf{F} = \pi M$ for some $M\in gr(R)$. Suppose that $j=0$. Since $\chi_1(M)$ holds, $\ul{\Ext}^i_R (R/R_{\geqslant n},M)$ is a finite $R$-module for each $i=1,2$ and together with \eqref{exactseq} it implies that $\omega\Sheaf{F} \cong \ul{H}^0(\Sheaf{F})$ is finite (recall $\omega$ from the paragraph after [\ref{notation}]). Now taking the $0$-graded part on both sides we get {\it (H4)} for $j=0$.

Suppose that $j\geqslant 1$. Since $R$ satisfies $\chi_{j+1}$, invoking proposition [\ref{chi_equiv}] we get

\beqn
\underset{\rightarrow}{\Lim} \ul{\Ext}^{j+1}_R (R/R_{\geqslant n},M)
\eeqn

is right bounded. Combining \eqref{cohom_desc} and \eqref{Ext_relation} this
equals $\ul{H}^j (\Sheaf{F})$. This immediately proves {\it (H5)} as
$\ul{H}^j(\Sheaf{F})_d = H^j(\Sheaf{F}[d])$. We now need left boundedness and local finiteness of $\ul{H}^j(\Sheaf{F})$ to finish the proof of {\it (H4)} for $j\geqslant 1$. These we have already observed (at least tacitly) but one can verify them by writing down a resolution of $R/R_{\geqslant n}$ involving finite sums of shifts of $R$, and then realising the cohomologies as subquotients of a complex of modules of the form $\ul{\Hom}_R (\underset{i=0}{\overset{p}{\oplus}} R[l_i],M)$.

\QED

Our discussion does not quite look complete unless we investigate the question of the ``dimension'' of the objects that we have defined.

\begin{center}
\underline{\bf Dimension of $\Proj\, R$}
\end{center}

\vspace{3mm}

The {\it cohomological dimension} of $\Proj\, R$, denoted by $cd(\Proj\, R)$, is defined to be

\begin{eqnarray*}
cd(\Proj\, R) :=
  \begin{cases}
     &\text{$sup\{i\; |\; {H}^i(\mathcal{M})\neq 0$ for some $\mathcal{M}\in \qgr(R)$\} if it is finite,} \\
     &\text{$\infty$ otherwise.}
  \end{cases}
\end{eqnarray*}

\begin{rem}
As ${H}^i$ commutes with direct limits one could have used $\QGr(R)$ in the definition of cohomological dimension.
\end{rem}

The following proposition gives us what we expect from a $\Proj$ construction regarding dimension and also provides a useful way of calculating it.

\begin{prop}

1. If $cd(\Proj\, R)$ is finite, then it is equal to $sup\{i\; |\; \underline{H}^i(R)\neq 0\}$.
\noindent

2. If the left global dimension of $R$ is $d<\infty$, then $cd(\Proj\, R)\leqslant d-1$.
\end{prop}

\noindent
{\bf Sketch of proof:}

1. Let $d$ be the cohomological dimension of $\Proj\, R$. It is obvious that $sup\{i\; |\; \underline{H}^i(R)\neq 0\}\leqslant d$. We need to prove the other inequality. So we choose an object for which the supremum is attained, {\it i.e.,} $\mathcal{M}\in \qgr(R)$ such that $H^d(\mathcal{M})\neq 0$ and, hence, $\underline{H}^d(\mathcal{M})\neq 0$. By the ampleness condition we may write down the following exact sequence:

\begin{eqnarray*}
0\map \mathcal{N}\map\underset{i=0}{\overset{p}{\oplus}} R[-l_i]\map \mathcal{M}\map 0
\end{eqnarray*} for some $\mathcal{N}\in \qgr(R)$. By the long exact sequence of derived functors $\underline{H}^i$ we have

\begin{eqnarray*}
\dots\map\underset{i=0}{\overset{p}{\oplus}}\underline{H}^d(R[-l_i])\map\underline{H}^d(\mathcal{M})\map\underline{H}^{d+1}(\mathcal{N})=0
\end{eqnarray*} This says that $\underline{H}^d(R[-l_i])\neq 0$ for some $i$ and hence, $\underline{H}^d(R)\neq 0$.

\vspace{1mm}
2. It has already been observed that $\underline{H}^i(M)\cong
\underset{n\rightarrow \infty}{\Lim}\underline{\Ext}^i(R_{\geqslant n},M)$ for
all $i\geqslant 0$. Now, if the left global dimension of $R$ is $d$, then
$\underline{\Ext}^j(N,M)=0$ for all $j>d$ and all $N,M\in \Gr(R)$. Putting
$N=R/R_{\geqslant n}$ and using \eqref{Ext_relation} we get $\underline{H}^d(M)=0$ for all $M\in \Gr(R)$. Therefore, $cd(\Proj\, R)\leqslant d-1$.

\QED

\begin{rem}
If $R$ is a Noetherian {\it AS-regular} graded algebra, then the {\it Gorenstein} condition can be used to prove that $cd(\Proj\, R)$ is actually equal to $d-1$.
\end{rem}

\begin{center}
\underline{\bf Some Examples (mostly borrowed from \cite{AZ})}
\end{center}

\vspace{3mm}

\begin{ex} {\it (Twisted graded rings)}

Let $\sigma$ be an automorphism of a graded algebra $A$. Then define a new multiplication $\ast$ on the underlying graded $k$-module $A = \underset{n}{\oplus} A_n$ by

\begin{eqnarray*}
a\ast b = a\sigma^n(b)
\end{eqnarray*} where $a$ and $b$ are homogeneous elements in $A$ and $\deg(a)=n$. Then algebra is called the {\it twist} of $A$ by $\sigma$ and it is denoted by $A^\sigma$. By \cite{ATV2} and \cite{Zha} $gr(A)\simeq gr(A^\sigma)$ and hence, $\proje(A)\simeq \proje(A^\sigma)$.

For example, if $A=k[x,y]$ where $\deg(x)=\deg(y)=1$, then any linear operator on the space $A_1$ defines an automorphism, and hence a twist of $A$. If $k$ is an algebraically closed field then, after a suitable linear change of variables, a twist can be brought into one of the forms $k_q[x,y]:=k\{x,y\}/(yx-qxy)$ for some $q\in k$, or $k_j[x,y]:=k\{x,y\}/(x^2+xy-yx)$. Hence, $\proje\, k[x,y]\simeq \proje\, k_q[x,y]\simeq \proje\, k_j[x,y]$. The projective scheme associated to any one of these algebras is the projective line $\proj^1$.
\end{ex}

\begin{ex} {\it (Changing the structure sheaf)}

Though the structure sheaf is a part of the definition of $\Proj$, one might
ask, given a $k$-linear abelian category $\mathcal{C}$, which objects
$\mathcal{A}$ could serve the purpose of the structure sheaf. In other words,
for which $\mathcal{A}$ do the conditions (H1), (H2) and (H3) of Theorem
[\ref{characterization}] hold? Since (H3) involves both the structure sheaf
and the polarization $s$, the answer may depend on $s$. We propose to
illustrate the possibilities by the simple example in which
$\mathcal{C}=\Mod(R)$ when $R=k_1\oplus k_2$, where $k_i = k$ for $i=1,2$ and
where $s$ is the automorphism which interchanges the two factors. The objects
of $\mathcal{C}$ have the form $V\simeq k_1^{n_1}\oplus k_2^{n_2}$, and the
only requirement for (H1), (H2) and (H3) is that both $r_1$ and $r_2$ not be zero simultaneously.

We have $s^n(V)=k_1^{r_2}\oplus k_2^{r_1}$ if $n$ is odd and $s^n(V)=V$ otherwise. Thus, if we set $\mathcal{A}=V$ and $A=\Gamma_h(\mathcal{C},\mathcal{A},s)$, then $A_n\simeq k_1^{r_1\times r_1}\oplus k_2^{r_2\times r_2}$ if $n$ is even, and $A_n\simeq k_1^{r_1\times r_2}\oplus k_2^{r_2\times r_1}$ otherwise. For example, if $r_1 =1$ and $r_2 =0$, then $A\simeq k[y]$, where $y$ is an element of degree $2$. Both of the integers $r_i$ would need to be positive if $s$ were the identity functor.
\end{ex}

\begin{ex} {\it (Commutative Noetherian algebras satisfy $\chi$ condition)}

Let $A$ be a commutative Noetherian $k$-algebra. Then the module structure on $\underline{\Ext}^n_A (A/A_+,M)$ can be obtained both from the right $A$-module structure of $A/A_+$ and that of $M$. Choose a free resolution of $A/A_+$, consisting of finitely generated free modules. The cohomology of this complex of finitely generated $A$-modules is given by the $\underline{\Ext}$'s, whence they are finite.

\end{ex}

\begin{ex} {\it (Noetherian AS-regular algebras satisfy $\chi$ condition)}

If $A$ is a Noetherian connected $\mathbb{N}$-graded algebra having global
dimension $1$, then $A$ is isomorphic to $k[x]$, where $\deg(x)=n$ for some
$n>0$, which satisfies the condition $\chi$ by virtue of the previous example. In higher dimensions we have the following proposition.

\begin{prop}
Let $A$ be a Noetherian AS-regular graded algebra of dimension $d\geqslant 2$ over a field $k$. Then $A$ satisfies the condition $\chi$.
\end{prop}

\noindent
{\bf Sketch of proof:}

$A$ is Noetherian and locally finite (due to finite $\GKdim$). For such an $A$
it is easy to check that $\underline{\Ext}^j(N,M)$ is a locally finite
$k$-module whenever $N,M$ are finite. $A_0$ is finite and hence,
$\underline{\Ext}^j(A_0 ,M)$ is locally finite for every finite $M$ and every
$j$. Since $A$ is connected graded, $A_0 = k$. For any $n$ and any finite
$A$-module $M$ we first show that $\underline{\Ext}^n(A_0,M) =
\underline{\Ext}^n(k,M)$ is bounded using induction on the projective
dimension of $M$. If $pd(M)=0$, then $M=\underset{i=0}{\overset{p}{\oplus}}
A[-l_i]$. By the {\it Gorenstein} condition (see definition
[\ref{Gorenstein}]) of an AS-regular algebra $A$, $\underline{\Ext}^n(k,A[-l_i])$ is bounded for each $i$. Therefore, so is $\underline{\Ext}^n(k,M)$. If $pd(M)>0$, we choose an exact sequence

\begin{eqnarray*}
0\map N\map P\map M\map 0
\end{eqnarray*} where $P$ is projective. Then $pd(N)=pd(M)-1$. By induction, $\underline{\Ext}^n(k,N)$ and $\underline{\Ext}^n(k,P)$ are bounded, hence, so is $\underline{\Ext}^n(k,M)$. Now $A/A_+$ is finite and we have just shown $\underline{\Ext}^n(k,M)$ is finite (since bounded together with locally finite implies finite); then $\Hom(A/A_+,\underline{\Ext}^n(k,M))\cong \underline{\Ext}^n(A/A_+,M)$ is locally finite and clearly bounded. Therefore, $\underline{\Ext}^n(A/A_+,M)$ is finite for every $n$ and every finite $M$.

\QED

\end{ex}

Most of these examples are taken directly from the original article by {\it Artin} and {\it Zhang} \cite{AZ}. There is a host of other examples on algebras satisfying $\chi$ up
to varying degrees, for which we refer the interested readers to
\cite{AZ}, \cite{Rog} and \cite{StZ1}. Also one should take a look
at \cite{ATV1} where these ideas, in some sense, germinated.
Finally, a comprehensive survey article by {\it J. T. Stafford} and {\it
M. van den Bergh} \cite{StV} should be consulted for further curiosities in the current state of affairs in \nc algebraic geometry.

%----------------------------Polishchuk's Work-------------------------------

\newpage

\section{Algebraic aspects of \nc tori}

\spc

The section is mostly based on the article {\it Noncommutative
two-tori with real multiplication as noncommutative projective
spaces} by {\it A. Polishchuk} \cite{Pol2}. Noncommutative two-tori with real multiplication will be
explained by {\it Jorge Plazas} and by now we know what we mean by
\nc projective varieties. In the following paragraph the gist of the
article has been provided (objects within quotes will be explained either by
{\it Jorge Plazas} or by {\it B. Noohi} or the reader is expected
to look it up for himself/herself). Interested readers
are also encouraged to take a look at the 4th chapter entitled {\it Fractional
  dimensions in homological algebra} of \cite{Man3}, where {\it Manin} gives a
very insightful overview of this work. We merely fill in some
details here for pedagogical reasons. 

One considers the category of ``holomorphic vector bundles'' on a
noncommutative torus $\mathbb{T}_\theta$, $\theta$ being a real parameter,
whose algebra of smooth functions is denoted $A_\theta$
\footnote{Usually $A_\theta$ is used to denote the algebra of
continuous functions and $\mathcal{A}_\theta$ is used to denote
the smooth ones. To ease LaTeX-ing, the algebra of smooth functions
has been consistently denoted by $A_\theta$.}. We always assume $\theta$ to be
irrational. In keeping with the general philosophy, $\mathbb{T}_\theta$ and
$A_\theta$ will be used interchangeably for a \nc torus. There is a
fully faithful functor \cite{PolSch} from the ``derived category of holomorphic
bundles on a noncommutative torus $A_\theta$'' to the derived category
of coherent sheaves on a complex elliptic curve $X$, denoted by $D^b(X)$ (refer
to the remark below \ref{derivedCategory}), sending holomorphic bundles to the
``heart'' of a ``t-structure'' depending on $\theta$, denoted by
$\mathcal{C}^\theta$. The elliptic curve $X$ is determined by the choice of a
``complex structure'' on $A_\theta$, depending on a complex parameter $\tau$ in the lower half plane {\it
i.e.,} $X = \CC/{(\ZZ + \ZZ\tau)}$. It follows from \cite{Pol1} that the
category of holomorphic vector bundles on $A_\theta$ is actually equivalent to
the heart $\mathcal{C}^\theta$ and the ``standard bundles'' end up being the so-called
``stable'' objects of $D^b(X)$. We also know that the heart has
``cohomological dimension'' $1$ and is derived equivalent to $D^b(X)$. The ``real multiplication'' of $A_\theta$ gives
rise to an auto-equivalence, say $F$, of $D^b(X)$, which preserves
the heart up to a shift. One knows when it actually preserves the heart, {\it viz.,} when the matrix inducing the real multiplication has positive real eigenvalues. Now by choosing a ``stable'' object, say $\Sheaf{G}$, in
$D^b(X)$, one can construct graded algebras from the triple $\mathcal{C}^\theta$, $\Sheaf{G}$ and
$F$ as described before (see \eqref{triple2alg}). Some criteria for
the graded algebras to be generated in degree $1$, quadratic and
``Koszul'' are also known. For the details one may refer to {\it e.g.,} \cite{Pla}.

\begin{rem}   \label{Con}
It is known that two noncommutative tori, say $A_\theta$ and $A_{\theta'}$, are
Morita equivalent if $\theta'=g\theta$ for some $g\in SL(2,\ZZ)$ \cite{Rie}.

The equivalence defined in \cite{PolSch} between the ``derived category of holomorphic bundles
on $A_\theta$'' and $D^b(X)$ actually sends the holomorphic bundles on
$A_\theta$ to $\mathcal{C}^{-\theta^{-1}}$ up to some shift and the
real multiplication on $A_\theta$ descends to an element $F\in
Aut(D^b(X))$, which preserves $\mathcal{C}^{-\theta^{-1}}$ (up to some shift). This is
not too bad, as $\begin{pmatrix}
                   0 & 1 \\
                   -1 & 0 \\
                 \end{pmatrix}\theta = -\theta^{-1}$ (action by fractional
linear transformation), which says that $A_{-\theta^{-1}}$ is Morita
equivalent to $A_\theta$. The generators of $SL(2,\ZZ)$ are $g=\begin{pmatrix}
                                    1 & 1 \\
                                    0 & 1 \\
                                   \end{pmatrix}$ and $h=\begin{pmatrix}
                                                           0 & 1 \\
                                                           -1 & 0 \\
                                                         \end{pmatrix}$. 
The  matrix $g$ acts by translation by $1$ and so up to Morita equivalence
                                                           $\theta$ may be
                                                           brought within the interval
                                                           $[0,1]$ and
                                                           $\theta$, being
                                                           irrational,
                                                           $\theta\in (0,1)$.

The image of this interval under $x\longmapsto -x^{-1}$ is $(-\infty,-1)$. We
label the noncommutative torus by $A_{-\theta^{-1}}$, $-\theta^{-1}\in(-\infty,-1)$, so that when we pass on to the heart
  $\mathcal{C}^{\theta}$, $\theta\in (0,1)$. 

\end{rem}

\spc

\begin{rem} \label{derivedCategory}
By the bounded derived category of coherent sheaves one should actually understand $D^b(\Coh(X))$. However, we may 
not be able to find injective resolutions in $\Coh(X)$. So the
precise category we want is $D^b(X):= D^b_{\Coh(X)}(\QCoh(X))$, {\it
i.e.}, the bounded derived category of complexes of quasi-coherent
sheaves on $X$ with the cohomology objects in $\Coh(X)$. According to 
{\it Lemma 2.3} of \cite{SeiTho} one knows that the two categories under consideration are equivalent. It is known that for a smooth curve $X$,
every object of $D^b(X)$ (which is a complex) is quasi-isomorphic to
the direct sum of its cohomologies. \footnote{By induction, it is
enough to show for complexes of length $2$. Let $F^{\bullet}\in
D^b(X)$.

 $F^{\bullet} = \dots\map 0\map F^{-1}\overset{f}{\map} F^{0}\map 0\map \dots$.

Consider the triangle: $ker\;
f[1]\overset{\theta}{\map}F^{\bullet}\map
cone\;\theta\overset{\xi}{\map} ker\; f[2]$. Check that in $D^b(X)$
$cone\; \theta$ is $coker\; f$. Now $\xi\in \Hom(coker\; f,ker\;
f[2]) = \Hom^{2} (coker\; f,ker\; f) = 0$. The last equality is due
to {\bf Fact 2} (it appears later on). So $F^{\bullet}[1] = coker\;
f[1]\oplus ker\; f[2]$. Note that $coker\; f$ and $ker\; f$ are the
cohomologies of $F^{\bullet}$.} This is not true in higher
dimensions.
\end{rem}

Let $X$ be an elliptic curve over the complex numbers. For $F\in \Coh(X)$, let $\rk(F)$ stand for the generic rank of $F$ and $\chi(F)$
for the Euler characteristic of $F$. Since $X$ has genus $1$, by
Riemann-Roch the degree of $F$ is the same as the Euler characteristic
of $F$, {i.e.,} $\deg(F)=\chi(F):= \dim_\CC \Hom(\sheaf_X,F)
-\dim_\CC \Ext^1(\sheaf_X,F)$. So the {\bf slope} of a coherent sheaf $F$, denoted by $\mu(F)$, which is just the
rational number $\frac{\deg(F)}{\rk(F)}$, also equals
$\frac{\chi(F)}{\rk(F)}$. The latter fraction is more suitable
for our purposes and, hence, we take that as the definition of slope. 
We define the rank (resp. the Euler characteristic) of a complex in the
derived category of $\Coh(X)$ as the
alternating sum of the ranks (resp. the Euler characteristics) of the individual
terms of the corresponding cohomology complex. Then the same definition as
above extends the notion of slope to the objects of the derived category. A coherent sheaf $F$ is called {\bf semistable} (resp. {\bf stable})
if for any nontrivial exact sequence $0\map
F'\map F\map F''\map 0$ one has
$\mu(F')\leqslant \mu(F)$
(resp. $\mu(F')<\mu(F)$) or equivalently
$\mu(F)\leqslant \mu(F'')$
(resp. $\mu(F)<\mu(F'')$).

It is well-known that every coherent sheaf on $X$ splits as a direct sum of
its torsion and torsion-free parts. Since $X$ is smooth, projective and of dimension $1$, every torsion-free
coherent sheaf is locally free and for any $F\in \Coh(X)$ there
exists a unique filtration \cite{HN}:

\begin{eqnarray} \label{HNF}
F=F_0\supset F_1\supset\dots\supset F_{n}\supset F_{n+1} = 0
\end{eqnarray} such that

\begin{itemize}

\item $F_i/F_{i+1}$ for $0\leqslant i\leqslant n$
  are semistable and

\item $\mu(F_0/F_1)<\mu(F_1/F_2)<\dots <\mu(F_n)$.

\end{itemize}

\noindent
The filtration above is called the {\bf Harder-Narasimhan filtration} of
$F$ and the graded quotients $F_i/F_{i+1}$ are called the {\bf semistable factors} of $F$. We set
$\mu_{min}(F) = \mu(F_0/F_1)$ and $\mu_{max}(F) = \mu(F_n)$.

One calls an object $F\in D^b(X)$ {\bf stable} if $F = V[n]$, where $V$ is
either a stable vector bundle (stable as above) or a coherent sheaf supported
at a point (the stalk is the residue field).

\subsection{t-structures on $D^b(X)$ depending on $\theta$}

One way to obtain t-structures is via ``torsion theories''. So let
us define a torsion pair $(\Coh_{>\theta},\Coh_{\leqslant\theta})$ in $\Coh(X)$.

\begin{eqnarray*}
\Coh_{>\theta} := \{ F\in \Coh(X)\; :\; \text{$\mu_{min}(F)>\theta$}\} \\
\Coh_{\leqslant\theta} := \{ F\in \Coh(X)\; :\; \text{$\mu_{max}(F)\leqslant\theta$}\}
\end{eqnarray*}

We consider the full subcategories generated by these
objects inside $\Coh(X)$. Notice that {\it torsion sheaves},
having slope = $\infty$, belong to $\Coh_{>\theta}$.

\spc

To show that this is indeed a torsion pair we need to verify two conditions.

\spc

{\it 1. $\Hom(T,F) = 0$ for all $T\in \Coh_{>\theta}$ and $F\in \Coh_{\leqslant\theta}$}

\begin{lem} \label{slopeCond}
Let $F$ and $F'$ be a pair of semistable bundles. Then $\mu (F)>\mu (F')$ implies that $\Hom(F,F')=0$.
\end{lem}

{\bf Proof:} Suppose $f:F\map F'$ is a nonzero morphism. Let $G$ be the image of $f$. Then $G$ is a quotient of $F$ and so one has $\mu (G)\geqslant\mu (F)$. On the other hand, $G$ is a torsion-free subsheaf of a vector bundle on a smooth curve and so it is locally free. Thus, one has $\mu (G)\leqslant \mu (F')$, which implies $\mu (F)\leqslant\mu (G)\leqslant\mu (F')$. Take the contrapositive to obtain the desired result. \QED

\spc

Let $T\in \Coh_{>\theta}$ and $F\in \Coh_{\leqslant\theta}$. Further, suppose $\sigma\in \Hom(T,F)$. Let us write down the {\it Harder-Narasimhan filtrations} of $T$ and $F$ respectively.

\begin{eqnarray*}
0\map T_0\map\cdots\map T_{m-1}\map T_m = T
\end{eqnarray*}

\beqn
0\map F_0\map\cdots\map F_{n-1}\map F_n = F
\eeqn

Restrict $\sigma : T\map F$ to $T_0$ and compose it with the canonical projection onto $F_n/F_{n-1}$. Now $T_0$ is a semistable factor of $T$ and $F_n/F_{n-1}$ that of $F$. Since $T\in \Coh_{>\theta}$ and $F\in \Coh_{\leqslant\theta}$, by the lemma above this map is $0$. So the image lies in $F_{n-1}$. Apply the same argument after replacing $F_n/F_{n-1}$ by $F_{n-1}/F_{n-2}$ to conclude that the image lies in $F_{n-2}$. Iterating this process we may conclude that $\sigma$ restricted to $T_0$ is $0$. So $\sigma$ factors through $T/T_0$. The {\it Harder-Narasimhan filtration} of $T/T_0$ is

\beqn
0\map T_1/T_0\map\cdots\map T_{m-1}/T_0\map T_m/T_0 = T/T_0
\eeqn

This filtration has the same semistable factors as that of $T$ and so they
satisfy the conditions of the {\it Harder-Narasimhan filtration}. So by the
uniqueness of the {\it Harder-Narasimhan filtration} this is that of $T/T_0$. Iterate the argument above after replacing $T$ by $T/T_0$ and taking the induced map of $\sigma$ between $T/T_0$ and $F$ to conclude that $\sigma$ vanishes on $T_1/T_0$. But $\sigma$ also vanishes on $T_0$. So it must vanish on $T_1$. Repeating this argument finitely many times one may show that $\sigma$ vanishes on the whole of $T$.

\spc

{\it 2. For every $F\in \Coh(X)$ there should be an exact sequence (necessarily
  unique up to isomorphism)

\begin{eqnarray*}
0\map t(F)\map F\map F/t(F)\map 0
\end{eqnarray*}

such that $t(F)\in \Coh_{>\theta}$ and $F/t(F)\in \Coh_{\leqslant\theta}$.}

\spc

{\bf Proof:} Let $0\subset F_1\subset\dots\subset F_{n-1}\subset F_n = F$ be the Harder-Narasimhan filtration of $F$. Let $i$ be the unique integer such that $\mu (F_i/F_{i-1})>\theta$ and $\mu (F_{i+1}/F_i)\leqslant\theta$. Then set $t(F) = F_i$. It is easy to see that $F_i\in \Coh_{>\theta}$ and $F/F_i\in \Coh_{\leqslant\theta}$. By the way, if no such $i$ exists, then $F$ is already either an element of $\Coh_{>\theta}$ or $\Coh_{\leqslant\theta}$. \QED

\spc

\begin{fact} \label{newTorsion} (see for instance \cite{HRS})
Let $(\mathcal{T},\mathcal{F})$ be a torsion pair on an abelian category $\mathcal{A}$. Let $\mathcal{C}$ be the heart of the associated t-structure. Then $\mathcal{C}$ is an abelian category, equipped with a torsion pair $(\mathcal{F}[1],\mathcal{T})$.
\end{fact}

Recall that the {\bf cohomological dimension} (perhaps, global dimension is a 
more appropriate term) of an abelian category $\mathcal{A}$ is the minimum integer $n$ such that $\Ext^i (A,B) = 0$ for all $A,B\in \mathcal{A}$ and for all $i>n$ and $\infty$ if no such $n$ exists.

\begin{fact} \label{Fact2} \cite{Ser}
If $X$ is a smooth projective curve ({\it i.e.}, $\dim\,X = 1$), then the cohomological dimension of $\Coh(X)$ is $1$.
\end{fact}

Now we shall associate a t-structure to this torsion pair (see for instance \cite{HRS}) as follows:

\begin{eqnarray*}
D^{\theta,\leqslant 0}\; :=\; \{K\in D^b(X)\; :\; H^{>0}(K) = 0, H^0 (K)\in \Coh_{>\theta}\} \\
D^{\theta,\geqslant 1}\; :=\; \{K\in D^b(X)\; :\; H^{<0}(K) = 0, H^0 (K) \in \Coh_{\leqslant\theta}\}
\end{eqnarray*}

\spc

It is customary to denote $D^{\theta,\leqslant 0}[-n]$ by $D^{\theta,\leqslant n}$ and $D^{\theta,\geqslant 0}[-n]$ by $D^{\theta,\geqslant n}$. Let $\mathcal{C}^\theta\; :=\; D^{\theta,\leqslant 0} \cap D^{\theta,\geqslant 0}$ be the heart of the t-structure, which is known to be an abelian category. An interesting thing is that $(\Coh_{\leqslant\theta}[1],\Coh_{>\theta})$ defines a torsion pair on $\mathcal{C}^\theta$ (refer to {\bf Fact 1} above). As a matter of convention, the family of t-structures is extended to $\theta = \infty$ by putting the standard t-structure on it, whose heart is just $\Coh(X)$.

\spc

Our next aim is to show that $\mathcal{C}^\theta$ has cohomological dimension $1$.

\spc
\noindent
{\bf ASIDE on Serre Duality:} Let $X$ be a smooth projective scheme of dimension $n$. Then there is a {\it dualizing sheaf} $\omega$ such that one has natural isomorphisms

\begin{eqnarray*}
H^i(X,F)\cong \Ext^{n-i}(F,\omega)^\ast
\end{eqnarray*}

where $F$ is any coherent sheaf on $X$.

\begin{rem}
The definition of a {\it dualizing sheaf} exists for all proper schemes. For nonsingular projective varieties it is known that the dualizing sheaf is isomorphic to the canonical sheaf.
\end{rem}

\begin{defn} \label{Serrefunctor} (\cite{BO1} Defn. 1.2., Prop. 1.3. and Prop. 1.4.)
Let $\mathcal{D}$ be a $k$-linear triangulated category with finite dimensional $\Hom$'s. An auto-equivalence $S:\mathcal{D}\map \mathcal{D}$ is called a {\bf Serre functor} if there are bi-functorial isomorphisms

\begin{eqnarray*}
\Hom_{\mathcal{D}} (A,B)\cong \Hom_{\mathcal{D}} (B,SA)^\ast
\end{eqnarray*}

which are natural for all $A,B\in \mathcal{D}$.
\end{defn}

{\it Bondal} and {\it Kapranov} have shown that in a reasonable manner Serre Duality of a smooth projective scheme $X$ can be reinterpreted as the existence of a Serre functor [if it exists it is unique up to a graded natural isomorphism] on $D^b (X)$ \cite{BK}.

It is also known that for a smooth projective variety of dimension $n$ the
Serre functor is ${}_{-}\otimes\omega_X[n]$. For an elliptic curve $X$ the
Serre functor will be just the translation functor $[1]$ (since $\dim\, X = 1$ and the canonical sheaf $\omega_X$ of an elliptic curve is trivial).

\spc

\begin{lem} \label{ssDecomp}
Any $F\in \Coh(X)$ is isomorphic to the direct sum of its semistable factors.
\end{lem}

\spc
\noindent
{\bf Proof:} 

The proof is by induction on the length of the
Harder-Narasimhan filtration and for simplicity we treat only the case of
length two. The category $\Coh(X)$ has the so-called {\it Calabi-Yau
property}, which says that $\Ext^1(F,G)\cong \Hom(G,F)^\ast$ for all $F,G\in \Coh(X)$ (this follows from
Serre Duality as discussed above). Let $F\in \Coh(X)$ and let 

\beqn
0\subset F_1\subset F_2=F
\eeqn be its Harder-Narasimhan filtration.

Then its semistable factors are $F_1$ and $F/F_1=:G$. Thus we obtain an exact
sequence 

\beqn
0\map F_1\map F\map G\map 0
\eeqn where $F_1$ and $G$ are semistable. By the properties of the Harder-Narasimhan filtration we have $\mu(G)<\mu(F_1)$. Due to the Calabi-Yau
property we know that $\Ext^1(G,F_1)\cong \Hom(F_1,G)^\ast$. From Lemma
\ref{slopeCond} it follows that $\Hom(F_1,G)=0$ and hence $\Ext^1(G,F_1)=0$. Therefore the short exact sequence above splits and we
obtain $F\cong F_1\oplus G$. \QED

\begin{prop}
$\mathcal{C}^\theta$ has cohomological dimension $1$.
\end{prop}

\spc
\noindent
{\bf Proof:}

First of all, observe that it is enough to show $\Hom^{>1}_{D^b(X)} (A,B) = 0$ for all $A,B$ belonging to $\Coh_{>\theta}$ and $\Coh_{\leqslant\theta}$ only. Now $\Coh(X)$ has cohomological dimension $1$ ({\bf Fact 2}) and $\Coh_{>\theta}$ and $\Coh_{\leqslant\theta}$ are full subcategories of $\Coh(X)$. So if $A,B$ were both either in $\Coh_{>\theta}$ or in $\Coh_{\leqslant\theta}[1]$ then there would have been nothing to prove. Let $A\in \Coh_{>\theta}$ and $B\in \Coh_{\leqslant\theta}[1]$.

Then $\Hom^i_{D^b(X)} (A,B) = \Hom^{i+1}_{D^b(X)} (A,B[-1])$. But $B\in \Coh_{\leqslant\theta}$, which is a subcategory of $\Coh(X)$. So $\Hom^{i+1}_{D^b(X)} (A,B[-1]) = 0$ for all $i\geqslant 1$. On the other hand,

\begin{eqnarray*}
\begin{split}
\Hom^i_{D^b(X)} (B,A) &= \Hom^{i-1}_{D^b(X)} (B[-1],A) \\
                     &= \Hom_{D^b(X)} (B[-1],A[1-i]) \\
                     &\cong \Hom_{D^b(X)} (A[1-i],B)^\ast \hspace{3mm} \text{(use Serre functor = $[1]$ as explained before)} \\
                     &= \Hom_{D^b(X)} (A[2-i],B[-1])^\ast \\
                     &= \Hom^{2-i}_{D^b(X)} (A,B[-1])^\ast
\end{split}
\end{eqnarray*}

\noindent So, for $i>2$, $\Hom^i_{D^b(X)} (B,A)$ is evidently $0$.
Due to the first axiom of a torsion pair, $\Hom_{\Coh(X)} (A,B[-1]) =
0$ when $i=2$. \hspace{53mm} \QED

\begin{prop}
The categories $\mathcal{C}^\theta$ and $\Coh(X)$ are derived equivalent, {\it
  i.e.,} $D^b(\mathcal{C}^\theta)\cong D^b(X)$. 
\end{prop}

\spc
\noindent
{\bf Proof:}
It is known that if a torsion pair $(\mathcal{T},\mathcal{F})$ in an abelian
category $\mathcal{A}$ is {\it cotilting}, {\it i.e., every object of
  $\mathcal{A}$ is a quotient of an object in $\mathcal{F}$}, then the heart
of the t-structure induced by the torsion pair is derived equivalent to
$\mathcal{A}$ (see for instance Proposition 5.4.3. and the remark thereafter
in \cite{BonVdb}). Thus, it is enough to check that the torsion pair
$(\Coh_{>\theta},\Coh_{\leqslant\theta})$ is cotilting. 

Given any $F\in \Coh(X)$ we need to produce an object in
$\Coh_{\leqslant\theta}$ which surjects onto $F$. Let $L$ be an ample line
bundle on $X$, {\it i.e.,} $\deg(L)>0$. By Serre's theorem one may twist $F$ by a large
enough power of $L$ such that it becomes generated by global sections, {\it
  i.e.,} the quotient of a free sheaf. In other
words, there exists $N>0$ large enough such that for all $n>N$ there is an
epimorphism ${\oplus}_{i\in I}\sheaf_X\map F\otimes L^n$, $I$ finite. One may twist it
back to obtain an epimorphism ${\oplus}_{i\in I}\check{L}^{n}\map F$, where
$\check{L}$ is the dual line bundle. This shows that there exists an epimorphism onto $F$ from a finite direct sum of copies of
$\check{L}^n$. Since $\deg(\check{L}^n)=-n.\deg(L)<0$ it is possible to make the
slope of $\check{L}^n$, which is equal to $\deg(\check{L}^n)$, less than
$\theta$ by choosing a large enough $n$. Being a line bundle $\check{L}^n$ is clearly semistable and we observe that the direct sum of copies of
$\check{L}^n$ lies in $\Coh_{\leqslant\theta}$.  \QED

\begin{rem}
In fact, all ``bounded'' t-structures on $D^b(X)$ come from some cotilting torsion
pair in $\Coh(X)$. All such cotilting torsion pairs (up to an action of
$Aut(D^b(X))$) have been listed in \cite{GKR}. I am especially thankful to
{\it S. A. Kuleshov} for explaining to me the above argument. 
\end{rem}

\spc

\noindent
Since $\mathcal{C}^\theta$ has cohomological dimension $1$, if it were
equivalent to $\Coh(Y)$ for some $Y$, then $Y$ had better be a smooth curve
({\it c.f.,} {\bf Fact 2}). The problem of dealing with categories of holomorphic bundles on $A_\theta$ has been reduced to studying t-structures on $D^b(X)$.

\spc

We have already seen some technical conditions involving a
categorical incarnation of ``ampleness'' \ref{ample} to verify when
a given $k$-linear ($k = \CC$ now) abelian category is of the
form $\Proj\, R$ for some graded $k$-algebra $R$. One of the
requirements of \ref{characterization} is that the category be
locally Noetherian ({\it i.e.,} the category has a Noetherian set of
generators). Unfortunately, this condition fails to be true in our
situation.

\begin{prop}
$\theta$ irrational implies that every nonzero object in $\mathcal{C}^\theta$ is not Noetherian.
\end{prop}

We would still like to say that what we have seen so far was not
entirely useless. Recall that in our discussion of $\Proj\, R$ after
{\it Artin} and {\it Zhang} we had assumed our graded algebra to be
right Noetherian. {\it Polishchuk} has shown that even if one
dispenses with the Noetherian assumption there is a way to recover
{\it Serre's theorem} \ref{characterization}. He gives an analogue
of an ``ample sequence of objects'' and proves that if a $k$-linear
abelian category has an ample sequence of objects then it is
equivalent to ``$\cohproj\, R$'', where $R$ is a ``coherent''
$\ZZ$-algebra. Unfortunately the words in quotes in the previous
sentence will not be explained anymore. Interested readers are
encouraged to look them up from \cite{Pol3}. Finally, as an apt
culmination of all our efforts we have the following theorem \cite{Pol2}.

\begin{thm}[Polishchuk]
For every quadratic irrationality $\theta\in\RR$ there exists an
auto-equivalence $F: D^b(X)\map D^b(X)$ preserving
$\mathcal{C}^\theta$ and a stable object $\Sheaf{G}\in
\mathcal{C}^\theta$ such that the sequence $(F^n\Sheaf{G},n\in\ZZ)$
is ample (in the modified sense of {\it Polishchuk}). Hence, the
corresponding algebra $A_{F,\Sheaf{G}}:= \Gamma_h
(\mathcal{C}^\theta,\Sheaf{G},F)$ is right coherent and
$\mathcal{C}^\theta\simeq \cohproj\, A_{F,\Sheaf{G}}$.
\end{thm}

\begin{rem}
There is some anomaly in the choice of the algebra $R$, whose
$\cohproj\, R$ should be equivalent to $\mathcal{C}^\theta$. However,
even in the commutative case one can show that if $S$ and $S'$ are
two graded commutative rings, such that $S_n \cong S'_n$ for all
$n>>0$, then $\Proj\, S \simeq \Proj\, S'$ (commutative $\Proj$
construction).
\end{rem}

\spc \noindent {\bf Final Remark:} ``Another perspective for the
future work is to try to connect our results with Manin's program in
\cite{Man3} to use noncommutative two-tori with real multiplication
for the explicit construction of the maximal Abelian extensions of
real quadratic fields''. --------- {\it A. Polishchuk}.

\bigskip

\noindent
{\bf Acknowledgements.} The author would like to acknowledge numerous useful
discussions with {\it Matilde Marcolli}, whose constant encouragement has also been pivotal to this work. The author would also like to
thank {\it Arend Bayer, I. Burban, Nikolai Durov, C. Kaiser, Behrang Noohi, Jorge
Plazas} and {\it Ramesh Sreekantan} for helpful discussions at various
stages. Finally, the author would like to express his gratitude to
{\it A. Greenspoon} for his careful proof-reading of the script, which helped remove several errors and imprecisions. 

%---------------bibliography-------------------------------------------

\bibliographystyle{alpha}
\bibliography{talk}

\end{document}